\documentclass[11pt]{article}
\usepackage{amsmath,amssymb,latexsym,epsf,mathabx}
\usepackage[all]{xy}
\usepackage{mathrsfs}
\usepackage{soul}
\usepackage{color}

\begin{document}

\newcommand{\nc}{\newcommand}
\newtheorem{Th}{Theorem}[section]
\newtheorem{Def}[Th]{Definition}
\newtheorem{Lem}[Th]{Lemma}
\newtheorem{Con}[Th]{Condition}
\newtheorem{Pro}[Th]{Proposition}
\newtheorem{Cor}[Th]{Corollary}
\newtheorem{Rem}[Th]{Remark}
\newtheorem{Exm}[Th]{Example}
\newtheorem{Sc}[Th]{}
\def\Pf#1{{\noindent\bf Proof}.\setcounter{equation}{0}}
\def\bskip#1{{ \vskip 20pt }\setcounter{equation}{0}}
\def\sskip#1{{ \vskip 5pt }\setcounter{equation}{0}}
\def\mskip#1{{ \vskip 10pt }\setcounter{equation}{0}}
\def\bg#1{\begin{#1}\setcounter{equation}{0}}
\def\ed#1{\end{#1}\setcounter{equation}{0}}

\soulregister\cite7 
\soulregister\citep7 
\soulregister\citet7 
\soulregister\ref7


\title{\bf  Quotients of extriangulated categories induced by selforthogonal subcategories
\thanks{Supported by the Natural Science Foundation of Universities of Anhui (Grant No.2023AH050950), the Top talent project of AHPU in 2020 (Grant No.S022021055),
the National Natural Science Foundation of China (Grant Nos.11801004, 12101003 and 12271249), the Natural Science Foundation of Anhui province (Grant No.2108085QA07) and
the Startup Foundation for Introducing Talent of AHPU (Grant No.2020YQQ067).}}
\smallskip
\author{ Peiyu Zhang, Yiwen Shi, Dajun Liu, Li Wang and Jiaqun Wei\\
\footnotesize ~E-mail:~zhangpy@ahpu.edu.cn, 1473008606@qq.com, ldjnnu2017004@163.com\\
\footnotesize ~wl04221995@163.com, weijiaqun5479@zjnu.edu.cn; \\
}

\date{}
\maketitle
\baselineskip 15pt
%
%
\begin{abstract}
\vskip 10pt%
Let $\mathcal{C}$ be an extriangulated category. We prove that two quotient categories of extriangulated categories
induced by selforthogonal subcategories are equivalent to module categories by
restriction of two functors $\mathbb{E}$ and $\mathrm{Hom}$, respectively. Moreover, if the selforthogonal subcategory is
contravariantly finite, then one of the two quotient categories is abelian. This result can be regarded as a generalization of Demonet-Liu and Zhou-Zhu.

\mskip\




\end{abstract}
%
\vskip 30pt

\section{Introduction}
Extriangulated categories were introduced by Nakaoka and Palu as a simultaneous generalization of exact categories and triangulated categories in \cite{NP}.
The class of extriangulated categories not only contains exact categories and extension closed subcategories of triangulated categories as examples, but it is also closed under taking some quotients.
There are many other examples for extriangulated categories which are neither exact categories nor triangulated categories, see \cite{NP} and \cite{ZZ}.
Many scholars have done a lot of research on extriangulated categories, for instance \cite{HZ,HZZ1,HZZ,LN,WWZ,ZB} and so on.

Let $\mathcal{T}$ be a cluster-tilting subcategory in a cluster
category $\mathcal{C}$. Buan-Marsh-Reiten \cite{BMR} proved that the quotient category $\mathcal{C}/\mathcal{T}[1]$ is
equivalent to the category of finitely presented modules $\mathrm{mod}\mathcal{T}$, where $\mathrm{mod}\mathcal{T}$ is abelian.
Let $\mathcal{C}$ be an exact category with enough projective objects and enough injective objects, and $\mathcal{T}$ a cluster-tilting subcategory of $\mathcal{C}$.
Demonet-Liu \cite{DL} proved that the quotient category $\mathcal{C}/\mathcal{T}$ is equivalent to the category of finitely presented
modules $\mathrm{mod}\mathcal{\underline{T}}$, where $\mathcal{\underline{T}}$ is the stable category of $\mathcal{T}$ by projective objects and $\mathrm{mod}\mathcal{\underline{T}}$ is abelian.
Let $\mathcal{C}$ be an extriangulated category with enough projective objects and enough injective objects, and $\mathcal{X}$ a rigid subcategory of $\mathcal{C}$.
Zhou-Zhu \cite{ZZ1} proved that the two quotient categories induced by $\mathcal{X}$ are equivalent to the category of finitely presented
modules $\mathrm{mod}\mathcal{\underline{X}}$ and $\mathrm{mod}\mathcal{\overline{X}}$, respectively,
where $\mathcal{\underline{X}}$ and $\mathcal{\overline{X}}$ are the stable category of $\mathcal{X}$ by projective objects and injective objects, respectively.
And $\mathrm{mod}\mathcal{\underline{X}}$ is abelian.

In this paper, let $\mathcal{C}$ be an extriangulated category with enough projective objects and enough injective objects, and $\mathcal{X}$ a selforthogonal subcategory of $\mathcal{C}$.
We denoted by $\underline{\mathcal{X}}$ the stable category which is quotient category of $\mathcal{X}$ by some strong contravariantly finite subcategory $\mathcal{U}$ of $\mathcal{X}$.
Similarly, for a strong covariantly finite subcategory $\mathcal{V}$ of $\mathcal{X}$, we denoted by $\overline{\mathcal{X}}=\mathcal{X}/[\mathcal{V}]$.
We mainly consider the two equivalences of categories: $\mathrm{mod}\mathcal{\underline{X}}\cong ?$ and $\mathrm{mod}\mathcal{\overline{X}}\cong ?$.
The conclusion is as follows.

\bg{Th}\label{Th0}
Let $\mathcal{C}$ be an extriangulated category with enough projective objects and enough injective objects, and $\mathcal{X}$ a selforthogonal subcategory of $\mathcal{C}$.

$(1)$ If $\mathcal{U}$ is a strong contravariant finite subcategory in $\mathcal{X}$, then $\mathcal{X}^{\mathcal{U}}_{L}/[\mathcal{X}]\cong \mathrm{mod}\mathcal{\underline{X}}$;

$(2)$ If $\mathcal{V}$ is a strong covariantly finite subcategory in $\mathcal{X}$, then $\mathcal{X}^{R}_{\mathcal{V}}/[\Omega^{-1}_{\mathcal{V}}\mathcal{X}]\cong \mathrm{mod}\mathcal{\overline{X}}$.
\ed{Th}



\section{Preliminaries}

First, we recall some definitions and relative properties of extriangulated category from \cite{NP}. Throughout this section, we assume that $\mathcal{C}$ is an additive category.
We denote by $\mathrm{Hom}_{\mathcal{C}}(A,B)$ or $\mathcal{C}(A,B)$ the set of morphisms from $A$ to $B$ in $\mathcal{C}$.
We use $\mathrm{Ab}$ to denote the category of abelian groups.

\bg{Def}\rm(\cite[Definition~2.1]{NP})\label{}
Assume that that $\mathcal{C}$ is equipped with an additive bifunctor $\mathbb{E}$: $\mathcal{C}^{op}\times \mathcal{C}\longrightarrow \mathrm{Ab}$, where $\mathrm{Ab}$ is the category of Abelian groups.
For any pair of objects $A,C\in\mathcal{C}$, an element $\delta \in \mathbb{E}(C,A)$ is
called an $\mathbb{E}$-extension. Thus formally, an $\mathbb{E}$-extension is a triplet $(A,\delta,C)$. If $\delta=0$, we called the split $\mathbb{E}$-extension.
\ed{Def}

For any morphism $a\in\mathcal{C}(A,A')$ and $c\in\mathcal{C}(C',C)$, we have $\mathbb{E}$-extensions
$\mathbb{E}(C,a)(\delta)\in\mathbb{E}(C,A')$ and $E(c,A)(\delta)\in\mathbb{E}(C',A)$.
We abbreviately denote them by $a_{\ast}\delta$ and $c^{\ast}\delta$ respectively.
A morphism $(a,c)$: $\delta\longrightarrow\delta'$ of $\mathbb{E}$-extensions is a pair of morphisms $a\in\mathcal{C}(A,A')$ and
$c\in\mathcal{C}(C,C')$ in $\mathcal{C}$ satisfying the equality $a_{\ast}\delta=c^{\ast}\delta'$.

\bg{Def}\rm(\cite[Definition~2.7]{NP})\label{}
Let $A,C\in \mathcal{C}$ be any pair of objects. Two sequences of morphisms $\xymatrix{A\ar[r]^{x} &B\ar[r]^{y} &C}$ and $\xymatrix{A\ar[r]^{x'} &B'\ar[r]^{y'} &C}$ are said to be
equivalent if there is an isomorphism $b\in\mathcal{C}(B,B')$ such that $bx=x'$ and $y'b=y$. We denote the equivalence class of $\xymatrix{A\ar[r]^{x} &B\ar[r]^{y} &C}$
by $[\xymatrix{A\ar[r]^{x} &B\ar[r]^{y} &C}]$.
\ed{Def}

%

\bg{Def}\rm(\cite[Definition~2.9]{NP})\label{}
Let $\mathfrak{s}$ be a correspondence which associates an equivalence
class $\mathfrak{s}(\delta)=[\xymatrix{A\ar[r]^{x} &B\ar[r]^{y} &C}]$ to any $\mathbb{E}$-extension $\delta\in\mathbb{E}(C,A)$.
This $\mathfrak{s}$ is called a realization of $\mathbb{E}$, if it satisfies the following condition $(\ast)$. In this
case, we say that sequence $[\xymatrix{A\ar[r]^{x} &B\ar[r]^{y} &C}]$ realizes $\delta$, whenever it satisfies
$\mathfrak{s}(\delta)=[\xymatrix{A\ar[r]^{x} &B\ar[r]^{y} &C}]$.
\ed{Def}

$(\ast)$ Let $\delta \in\mathbb{E}(C,A)$ and $\delta'\in\mathbb{E}(C',A')$ be any pair of $\mathbb{E}$-extensions, with
$\mathfrak{s}(\delta)=[\xymatrix@C=0.5cm{A\ar[r]^{x} &B\ar[r]^{y} &C}]$ and $\mathfrak{s}(\delta')=[\xymatrix{A'\ar[r]^{x'} &B'\ar[r]^{y'} &C'}]$.
Then, for any morphism $(a,c)$: $\delta\longrightarrow \delta'$, there exists $b\in\mathcal{C}(B,B')$ such that $bx=x'a$ and $y'b=cy$.

\bg{Def}\rm(\cite[Definition~2.10]{NP})\label{}
A realization $\mathfrak{s}$ of $\mathbb{E}$ is said to be additive, if it satisfies the following conditions.

$(1)$ For any $A,C\in \mathcal{C}$, the split $\mathbb{E}$-extension $0\in\mathbb{E}(C,A)$ satisfies $\mathfrak{s}(0) = 0$.

$(2)$ For any pair of $\mathbb{E}$-extensions $\delta$ and $\delta'$,
$\mathfrak{s}(\delta\oplus\delta')=\mathfrak{s}(\delta)\oplus\mathfrak{s}(\delta')$ holds.
\ed{Def}

\bg{Def}(\cite[Definition~2.12]{NP})\label{etri}
We call the triplet $(\mathcal{C},\mathbb{E},\mathfrak{s})$ (simply, $\mathcal{ C}$) an extriangulated category if it satisfies the following conditions:

$\mathrm{(ET1)}$ $\mathbb{E}$: $C^{op}\times C\longrightarrow \mathrm{Ab}$ is an additive functor.

$\mathrm{(ET2)}$ $\mathfrak{s}$ is an additive realization of $\mathbb{E}$.

$\mathrm{(ET3)}$ Let $\delta\in\mathbb{E}(C,A)$ and $\delta'\in\mathbb{E}(C',A')$ be any pair of $\mathbb{E}$-extensions,
realized as $\mathfrak{s}(\delta)=[\xymatrix{A\ar[r]^{x} &B\ar[r]^{y} &C}]$ and $\mathfrak{s}(\delta')=[\xymatrix{A'\ar[r]^{x'} &B'\ar[r]^{y'} &C'}]$.
For any commutative square in $\mathcal{C}$
$$\xymatrix{
A\ar[r]^{x}\ar[d]^{a}&B\ar[r]^{y}\ar[d]^{b}&C\\
A'\ar[r]^{x'} &B'\ar[r]^{y'} &C'
}$$
there exists a morphism $(a,c)$: $\delta\longrightarrow\delta'$ which is realized by $(a,b,c)$.

$(\mathrm{ET3})$$^{\mathrm{op}}$  Dual of $\mathrm{(ET3)}$.

$(\mathrm{ET4})$ Let $(A,\delta,D)$ and $(B,\delta',F)$ be $\mathbb{E}$-extensions respectively realized by $\xymatrix{A\ar[r]^{f} &B\ar[r]^{f'} &D}$
and $\xymatrix{B\ar[r]^{g} &C\ar[r]^{g'} &F}$. Then there exist an object $E\in\mathcal{C}$, a commutative diagram
$$\xymatrix{
A\ar[r]^{f}\ar@{=}[d] &B\ar[r]^{f'}\ar[d]^{g}&D\ar[d]^{d}\\
A\ar[r]^{h} &C\ar[r]^{h'}\ar[d]^{g'} &E\ar[d]^{e}\\
&F\ar@{=}[r]&F
}$$
in $\mathcal{C}$, and an $\mathbb{E}$-extension $\delta''\in\mathbb{E}(E,A)$ realized by $\xymatrix{A\ar[r]^{h} &C\ar[r]^{h'} &E}$,
which satisfy the following compatibilities.

$(i)$ $\xymatrix{D\ar[r]^{d} &E\ar[r]^{e} &F}$ realizes $\mathbb{E}(F,f')(\delta')$,

$(ii)$ $\mathbb{E}(d,A)(\delta'') = \delta$,

$(iii)$ $\mathbb{E}(E,f)(\delta'') = \mathbb{E}(e,B)(\delta')$.

By (iii), $(f,e)$: $\delta''\longrightarrow\delta$ is a morphism of $\mathbb{E}$-extensions, realized by
$$(f,id_{C} ,e):~[\xymatrix{A\ar[r]^{h} &C\ar[r]^{h'} &E}]\longrightarrow [\xymatrix{B\ar[r]^{g} &C\ar[r]^{g'} &F}].$$
$\mathrm{(ET4)}$$^{\mathrm{op}}$ Dual of $\mathrm{(ET4)}$.
\ed{Def}

%
%

In Definition \ref{etri}, if the sequence $\xymatrix{A\ar[r]^{x} &B\ar[r]^{y} &C}$ realizes $\delta\in\mathbb{E}(C,A)$,
we call the pair $(\xymatrix{A\ar[r]^{x} &B\ar[r]^{y} &C},\delta)$ an $\mathbb{E}$-triangle, and write it in the following way,
$$\xymatrix{A\ar[r]^{x} &B\ar[r]^{y} &C\ar@{.>}[r]^{\delta}&}.$$
And the sequence $\xymatrix{A\ar[r]^{x} &B\ar[r]^{y} &C}$ is called a conflation, $x$ is called an inflation and $y$ is called
a deflation.


%
%
%
%

\bg{Def}$\rm(\cite[Definition~4.1]{HZZ})$\label{}
Let $\mathcal{C}$ be an extriangulated category. An object $P\in \mathcal{C}$ is called projective if for any  $\mathbb{E}$-triangle
$\xymatrix@C=0.5cm{A\ar[r] &B\ar[r] &C\ar@{.>}[r]^{\delta}&}$, the induced sequence of of abelian group
$0\longrightarrow \mathrm{Hom}_{\mathcal{C}}(P,~A)\longrightarrow\mathrm{Hom}_{\mathcal{C}}(P,~B)\longrightarrow\mathrm{Hom}_{\mathcal{C}}(P,~C)\longrightarrow0$
is exact in $\mathrm{Ab}$. Dually, we can define injective objects.
\ed{Def}

We denote $\mathcal{P}$ (resp. $\mathcal{I}$) the class of projective (resp. injective) objects of $\mathcal{C}$. An extrianglated category $\mathcal{C}$ is said to have enough projective (resp. injective) objects if for any object $A\in\mathcal{C}$ there is an $\mathbb{E}$-triangle
$\xymatrix{B\ar[r] &P\ar[r] &A\ar@{.>}[r]&}$
(resp. $\xymatrix{A\ar[r] &I\ar[r] &B\ar@{.>}[r]&}$)
with $P\in \mathcal{P}$ (resp. $I\in \mathcal{I}$).

\bg{Def}$\rm(\cite{HZZ})\label{}$
Let $\mathcal{W}$ be a class of objects in $\mathcal{C}$. An $\mathbb{E}$-triangle $\xymatrix@C=0.5cm{A\ar[r] &B\ar[r]&C\ar@{.>}[r]&}$ is
said to be $\mathcal{W}$-exact if for any $W\in \mathcal{W}$, the induced sequence of abelian group
$0\longrightarrow\mathcal{C}(C,~W)\longrightarrow\mathcal{C}(B,~W)\longrightarrow\mathcal{C}(A,~W)\longrightarrow0$
is exact. Dually we can define $\mathcal{W}$-coexact sequences.
\ed{Def}

%
%
%
Assume that $\mathcal{C}$ is an extriangulated category. By Yoneda's lemma, any $\mathbb{E}$-extension
$\delta\in \mathbb{E}(C, A)$ induces natural transformations \cite[Definition~3.1]{NP}
$$\delta_{\sharp}: \mathcal{C}(-,C)\Rightarrow\mathbb{E}(-,A) ~\mathrm{and} ~
\delta^{\sharp}: \mathcal{C}(A,-)\Rightarrow\mathbb{E}(C,-)$$
For any $X\in \mathcal{C}$, these $\delta_{\sharp}$ and $\delta^{\sharp}$ given as follows:

(1) $(\delta_{\sharp})_{X}$: $\mathcal{C}(X,C)\rightarrow\mathbb{E}(X,A)$; $f\longmapsto f^{\ast}\delta$.

(2) $(\delta^{\sharp})_{X}$: $\mathcal{C}(A,X)\rightarrow\mathbb{E}(C,X)$; $g\longmapsto g_{\ast}\delta$.

\bg{Lem}\label{}
Let $\mathcal{C}$ be an extriangulated category and
$$\xymatrix{A\ar[r]^{a} &B\ar[r]^{b} &C\ar@{.>}[r]^{\delta}&}$$
an $\mathbb{E}$-triangle. Then we have the following long exact sequence:
$$\xymatrix{\mathcal{C}(-,A)\ar[r]^{\mathcal{C}(-,a)}&\mathcal{C}(-,B)\ar[r]^{\mathcal{C}(-,b)}&
\mathcal{C}(-,C)\ar[r]^{\delta^{\sharp}_{-}}&\mathbb{E}(-,A)\ar[r]^{\mathbb{E}(-,a)}&\mathbb{E}(-,B)\ar[r]^{\mathbb{E}(-,b)}&\mathbb{E}(-,C),}$$
$$\xymatrix{\mathcal{C}(C,-)\ar[r]^{\mathcal{C}(b,-)}&\mathcal{C}(B,-)\ar[r]^{\mathcal{C}(a,-)}&
\mathcal{C}(A,-)\ar[r]^{\delta_{\sharp}^{-}}&\mathbb{E}(C,-)\ar[r]^{\mathbb{E}(b,-)}&\mathbb{E}(B,-)\ar[r]^{\mathbb{E}(a,-)}&\mathbb{E}(A,-).}$$
\ed{Lem}

\Pf. It follows from Proposition 3.3 and Proposition 3.11 in \cite{NP}.
\ \hfill $\Box$

\vskip 10pt
In the rest of this paper, $\mathcal{C}$ is also extriangulated category with enough projective objects and injective objects.
All subcategories of $\mathcal{C}$ are full subcategories and closed under isomorphisms and finite direct sums.

A subcategory $\mathcal{X}$ of $\mathcal{C}$ is said to be contravariantly finite in $\mathcal{C}$ if for any object $C$ in $\mathcal{C}$, there is some object $X\in \mathcal{X}$ and a morphism $f$: $X\to C$ such that for every $X'\in \mathcal{X}$ the following diagram is commutative:
$$\xymatrix{
&X'\ar[d]^{\forall}\ar@{.>}[dl]_{\exists}\\
X\ar[r]^{f}&C.\\
}$$
In this case, $f$ is called a right $\mathcal{X}$-approximation of $C$. Dually, we can define covariantly finite subcategories in $\mathcal{C}$ and left $\mathcal{X}$-approximations. Furthermore, a subcategory of $\mathcal{C}$ is said to be functorially finite in $\mathcal{C}$ if it is both contravariantly finite and covariantly finite in $\mathcal{C}$.

\bg{Def}\cite[Definition ~3.21]{ZZ}
Let $\mathcal{C}$ be an extriangulated category. A subcategory $\mathcal{X}$ of $\mathcal{C}$ is called strongly contravariantly finite, if for any object $C\in\mathcal{C}$, there
exists an $\mathbb{E}$-triangle
$$\xymatrix{C_{1}\ar[r] &X\ar[r]^{f}&C\ar@{.>}[r]&}$$
where $f$ is a right $\mathcal{X}$-approximation of $C$.

Dually, a subcategory $\mathcal{X}$ of $\mathcal{C}$ is called strongly covariantly finite, if for any object $C\in\mathcal{C }$, there exists an $\mathbb{E}$-triangle
$$\xymatrix{C\ar[r] &X\ar[r]^{g}&C^{1}\ar@{.>}[r]&}$$
where $f$ is a left $\mathcal{X}$-approximation of $C$.

A strongly contravariantly finite and strongly covariantly finite subcategory is called
strongly functorially finite.
\ed{Def}

\bg{Rem}\cite[Remark~2.9]{ZZ1}
Let $\mathcal{C}$ be an extriangulated category and $\mathcal{X}$ be a subcategory of $\mathcal{C}$.

$(1)$ The subcategory $\mathcal{X}$ is strongly contravariantly finite in $\mathcal{C}$
if and only if $\mathcal{X}$ is contravariantly finite in $\mathcal{C}$ containing $\mathcal{P}$.

$(2)$ The subcategory $\mathcal{X}$ is strongly covariantly finite in $\mathcal{C}$
if and only if $\mathcal{X}$ is covariantly finite in $\mathcal{C}$ containing $\mathcal{I}$.

$(3)$ The subcategory $\mathcal{X}$ is strongly functorially finite in $\mathcal{C}$
if and only if $\mathcal{X}$ is functorially finite in $\mathcal{C}$ containing $\mathcal{P}$ and $\mathcal{I}$.
\ed{Rem}

Let $\mathcal{C}$ be an extriangulated category and $f$: $A\to B$ a morphism in $\mathcal{C}$. A pseudokernel of $f$ is a morphism $g$: $K\to A$ such that for any $C\in \mathcal{C}$ the sequence of abelian groups
$\xymatrix{\mathcal{C}(C,K)\ar[r] &\mathcal{C}(C,A)\ar[r]&\mathcal{C}(C,B)}$
is exact. Equivalently, $g$ is a pseudokernel of $f$ if $fg = 0$ and for each morphism $h$: $C\to A$
such that $fh= 0$, there exists a (not necessarily unique) morphism $i$: $C\to K$ such that
$h=gi$. i.e., the following diagram is commutative.
$$\xymatrix{
&C\ar@{.>}[ld]_{i}\ar[rd]^{0}\ar[d]^{h}&\\
K\ar[r]^{g}&A\ar[r]^{f}&B.\\
}$$
Clearly, a pseudokernel $g$ of $f$ is a kernel if and only if $g$ is a monomorphism.

A $\mathcal{C}$-module is a contravariant functor $\mathbb{G}$: $\mathcal{C}\to \mathrm{Ab}$. Then $\mathcal{C}$-modules form an abelian category $\mathrm{Mod}\mathcal{C}$.
By Yoneda's lemma, representable functors are projective objects in $\mathrm{Mod}\mathcal{C}$.
We call $M\in\mathrm{Mod}\mathcal{C}$ coherent \cite{Au} if there exists an exact sequence
$$\xymatrix{\mathcal{C}(-,X_{1})\ar[r]&\mathcal{C}(-,X_{0})\ar[r]&M\ar[r]&0}$$
in $\mathrm{Mod}\mathcal{C}$ with $X_{0}$, $X_{1}\in\mathcal{C}$. We denote by $\mathrm{mod}\mathcal{C}$ the full subcategory of $\mathrm{Mod}\mathcal{C}$
consisting of coherent $\mathcal{C}$-modules. It is easily checked that $\mathrm{mod}\mathcal{C}$ is closed under cokernels
and extensions in $\mathrm{Mod}\mathcal{C}$. Moreover, $\mathrm{mod}\mathcal{C}$ is closed under kernels in $\mathrm{Mod}\mathcal{C}$ if and only
if $\mathcal{C}$ has pseudokernels. In this case, $\mathrm{mod}\mathcal{C}$ forms an abelian category (see \cite{Au}).

\section{Quotients of extriangulated categories I}
In this section, we mainly present the first half of Theorem \ref{Th0}. Firstly, we recall some basic definitions.
Let $\mathcal{C}$ be an extriangulated category with enough projective objects and injective objects,
and $\mathcal{X}$ a subcategory of $\mathcal{C}$. We denote by $[\mathcal{X}](A, B)$ the subgroup of
$\mathrm{Hom}_{\mathcal{C}}(A, B)$ consisting of morphisms which factors through objects in a subcategory $\mathcal{X}$ of $\mathcal{C}$.
The quotient category $\mathcal{C}/[\mathcal{X}]$ of $\mathcal{C}$ by a subcategory $\mathcal{X}$ is the category with the same objects
as $\mathcal{C}$ and the set of morphisms from $A$ to $B$
is the quotient of group of morphisms from $A$ to $B$ in $\mathcal{C}$ by the subgroup consisting of
morphisms factors through objects in $\mathcal{X}$

Assume that $\mathcal{U}$ is a subcategory of $\mathcal{X}$. If $\mathcal{U}$ is strong contravariantly finite
in $\mathcal{X}$, the stable category $\mathcal{\underline{X}}$ of $\mathcal{X}$ is the quotient category $\mathcal{X}/[\mathcal{U}]$,
i.e., the category which has same objects with $\mathcal{X}$ and morphisms are defined as follows
$$\mathrm{Hom}_{\mathcal{\underline{X}}}(A,B)=\mathrm{Hom}_{\mathcal{X}}(A,B)/[\mathcal{U}](A,B).$$
We denote by $\underline{f}$ the residue class in $\mathrm{Hom}_{\mathcal{\underline{X}}}(A,B)$ of a morphism $f$ of $\mathrm{Hom}_{\mathcal{X}}(A,B)$.

\bg{Pro}\label{prop1}
Let $\mathcal{C}$ be an extriangulated category and $\mathcal{U}\subseteq\mathcal{X}$ two subcategories of $\mathcal{C}$.
Assume that $\mathcal{U}$ is strong contravariantly finite in $\mathcal{X}$ and $\mathcal{X}$ is strong contravariantly finite in $\mathcal{C}$.
Then $\mathrm{mod}\mathcal{\underline{X}}$ is an abelian category.
\ed{Pro}

\Pf. It suffices to prove that $\mathcal{\underline{X}}$ has pseudokernels (see \cite[Section~2]{Au} or \cite[Chapter~III, Section~2]{Au1}).
For any morphism $f\in\mathcal{C}(A,B)$ with $A$, $B\in \mathcal{X}$. Since $\mathcal{U}$ is strong contravariantly finite in $\mathcal{X}$,
there is an $\mathbb{E}$-triangle
$$\xymatrix{B_{1}\ar[r]^{a} &U_{0}\ar[r]^{b}&B\ar@{.>}[r]&}$$
with $b$ a right $\mathcal{U}$-approximation of $B$.
By \cite[Corollary~3.16]{NP}, there is an $\mathbb{E}$-triangle
$$\xymatrix{B_{2}\ar[r]^{\binom{g}{h}~~~} &A\bigoplus U_{0}\ar[r]^{~~~(f,b)}&B\ar@{.>}[r]&}.$$
Since $\mathcal{X}$ is strong contravariantly finite in $\mathcal{C}$, there is a right $\mathcal{X}$-approximation $x$: $X\to B_{2}$.
Next, we will prove that $\underline{gx}$ is a pseudokernels of $\underline{f}$.
If $\underline{fx'}=0$, for a morphism $x'\in\mathrm{Hom}_{\mathcal{\underline{X}}}(X',A)$ with $X'\in \mathcal{X}$,
then there exists the following commutative diagram:
$$\xymatrix{
A\ar[r]^{f}&B\\
X'\ar[u]^{x'}\ar[r]^{d}&U\ar[u]^{c}
}$$
where $U\in \mathcal{U}$. Thus we have the following commutative diagram:
$$\xymatrix{
&&U\ar@{.>}[dl]_{e}\ar[d]^{c}&\\
B_{1}\ar[r]&U_{0}\ar[r]^{b}&B\ar@{.>}[r]&
}$$
since $b$ is a right $\mathcal{U}$-approximation of $B$.
Notice that $(f,b)\binom{x'}{-ed}=fx'-bed=fx'-cd=0$, then we have the following commutative diagram:
$$\xymatrix{
X'\ar@{.>}[d]^{y}\ar[dr]^{\binom{x'}{-ed}}&&&\\
B_{2}\ar[r]^{\binom{g}{h}~~~}&A\bigoplus U_{0}\ar[r]^{~~~(f,b)}&B\ar@{.>}[r]&
}$$
Since $x$ is a right $\mathcal{X}$-approximation of $B_{2}$, we have the following commutative diagram:
$$\xymatrix{
&X'\ar[d]^{y}\ar@{.>}[dl]_{z}\\
X\ar[r]^{x}&B_{2}
}$$
From above diagrams, we have that $x'=gy=gxz$. i.e., we have the following commutative diagram:
$$\xymatrix{
&X'\ar[d]^{\underline{x'}}\ar@{.>}[dl]_{\underline{z}}\ar[dr]^{0}&\\
X\ar[r]^{\underline{gx}}&A\ar[r]^{\underline{f}}&B.
}$$
Consequently, $\underline{gx}$ is a pseudokernels of $\underline{f}$.
\ \hfill $\Box$

\vskip 10pt
Let $\mathcal{X}$ be a selforthogonal subcategory of $\mathcal{C}$ (i.e., $\mathbb{E}(\mathcal{X},\mathcal{X})=0$)
and $\mathcal{U}$ be be a strong contravariant finite subcategory of $\mathcal{X}$.
We denote by $\mathcal{X}^{\mathcal{U}}_{L}$ the subcategory of objects $A$ of $\mathcal{C}$ which admits an $\mathbb{E}$-triangle:
$\xymatrix{A\ar[r] &X_{0}\ar[r]&X_{1}\ar@{.>}[r]&}$ with $X_{0}$ and $X_{1}$ in $\mathcal{X}$
and it is $\mathcal{U}$-exact. We consider the restriction of functor $\mathbb{E}(-, -)$ to $\mathcal{X}^{\mathcal{U}}_{L}$, denoted by $\mathbb{H}$:
$$\mathcal{X}^{\mathcal{U}}_{L}\longrightarrow \mathrm{Mod}\mathcal{\underline{X}}~; A\longmapsto \mathbb{H}(A)=\mathbb{E}(-,A)|_{\mathcal{X}}.$$
It is easy to see that $\mathcal{X}\subseteq\mathcal{X}^{\mathcal{U}}_{L}$, thus we have the quotient category $\mathcal{X}^{\mathcal{U}}_{L}/[\mathcal{X}]$
and denote by $[f]$ the residue class in $\mathcal{X}^{\mathcal{U}}_{L}/[\mathcal{X}]$ of a morphism $f$ of $\mathcal{X}^{\mathcal{U}}_{L}$.

Let $\pi$: $\mathcal{X}^{\mathcal{U}}_{L}\longrightarrow \mathcal{X}^{\mathcal{U}}_{L}/[\mathcal{X}]$ be projective functor. Since $\mathcal{X}$ is selforthogonal,
$\mathbb{H}(M)=\mathbb{E}(-,M)|_{\mathcal{X}}=0$. By the universal property of $\pi$, there is a functor
$$\mathbb{F}: \mathcal{X}^{\mathcal{U}}_{L}/[\mathcal{X}]\longrightarrow \mathrm{Mod}\mathcal{\underline{X}}$$
such that the following diagram commutes:
$$\xymatrix{
\mathcal{X}^{\mathcal{U}}_{L}\ar[rr]^{\pi}\ar[dr]_{\mathbb{H}}&&\mathcal{X}^{\mathcal{U}}_{L}/[\mathcal{X}]\ar@{.>}[dl]^{\mathbb{F}}\\
&\mathrm{Mod}\mathcal{\underline{X}}&.
}$$

The following result plays an important role in the proof of Theorem \ref{Th0}.

\bg{Lem}\label{lemma1}
Let $\mathcal{X}$ be a selforthogonal subcategory of an extriangulated category $\mathcal{C}$ and $\mathcal{U}$ a strong contravariant finite subcategory in $\mathcal{X}$.
For any $\mathcal{U}$-exact $\mathbb{E}$-triangle
$$(\dag):\xymatrix{A\ar[r]^{x} &X_{0}\ar[r]^{y}&X_{1}\ar@{.>}[r]&}$$
with $X_{0}$, $X_{1}\in \mathcal{X}$, then there is an exact sequence in $\mathrm{Mod}\mathcal{\underline{X}}$
$$\xymatrix{\mathrm{Hom}_{\mathcal{\underline{X}}}(-,X_{0})\ar[r]&\mathrm{Hom}_{\mathcal{\underline{X}}}(-,X_{1})\ar[r]&\mathbb{H}(A)\ar[r]&0.}$$
Thus we have that $\mathbb{H}(A)=\mathbb{F}(A)\in \mathrm{mod}\mathcal{\underline{X}}$.
\ed{Lem}

\Pf. For any $X\in \mathcal{X}$, applying the functor $\mathrm{Hom}_{\mathcal{C}}(X,-)$ to $\mathbb{E}$-triangle $(\dag)$, we can get the following
long exact sequence:
$$\xymatrix{ \mathrm{Hom}_{\mathcal{C}}(X,X_{0})\ar[r]^{\mathcal{C}(X,y)}& \mathrm{Hom}_{\mathcal{C}}(X,X_{1})\ar[r]^{~~~\delta}& \mathbb{E}(X,A)\ar[r]&\mathbb{E}(X,X_{0})=0}$$
since $\mathcal{X}$ is selforthogonal. Thus we have the following exact sequence in $\mathrm{Mod}\mathcal{X}$:
$$\xymatrix{ \mathrm{Hom}_{\mathcal{X}}(-,X_{0})\ar[r]& \mathrm{Hom}_{\mathcal{X}}(-,X_{1})\ar[r]& \mathbb{E}(-,A)\ar[r]&0}$$
For any $f\in[\mathcal{U}](X,X_{1})$, it is easy to prove that $f$ factors through $y$ since $\mathbb{E}$-triangle $(\dag)$ is $\mathcal{U}$-exact.
i.e., $f\in \mathrm{Im}~\mathcal{C}(X,y)$. Thus $\delta(f)=0$. So $\delta$ induce a morphism $\underline{\delta}$:
$$\mathrm{Hom}_{\mathcal{\underline{X}}}(X,X_{1})\longrightarrow \mathbb{E}(X,A)=\mathbb{H}(A)(X): \underline{f}\longmapsto \underline{\delta}(\underline{f})=\delta(f).$$
Based on the above discussion, the induce morphism $\underline{\delta}$ is well-define.
It is not difficult to see that $\underline{\delta}$ is surjective for any $X\in\mathcal{X}$.

Since $\delta\circ\mathrm{Hom}_{\mathcal{X}}(X,y)=0$, we can have that $\underline{\delta}\circ\mathrm{Hom}_{\mathcal{\underline{X}}}(X,\underline{y})=0$.
If $\underline{\delta}(\underline{g})=0$ for any $g\in \mathrm{Hom}_{C}(X,X_{1})$, then $\delta(g)=\underline{\delta}(\underline{g})=0$.
Thus there is a morphism $h\in\mathrm{Hom}_{\mathcal{C}}(X,X_{0})$ such that $g=yh$, and then $\underline{g}=\underline{y}\underline{h}$.
i.e., $\mathrm{Im}\mathrm{Hom}_{\mathcal{\underline{X}}}(X,\underline{y})=\mathrm{Ker} \underline{\delta}$. Consequently, we have the following exact
sequence in $\mathrm{mod}\underline{\mathcal{X}}$:
$$\xymatrix{ \mathrm{Hom}_{\mathcal{\underline{X}}}(-,X_{0})\ar[r]
& \mathrm{Hom}_{\mathcal{\underline{X}}}(-,X_{1})\ar[r]& \mathbb{E}(-,A)=\mathbb{H}(A)\ar[r]&0.}$$
So we complete this proof.
\ \hfill $\Box$

\vskip 10pt

Next, we will give the proof of the first half of Theorem 1.1.

\bg{Th}\label{Th1}
Let $\mathcal{X}$ be a selforthogonal subcategory of an extriangulated category $\mathcal{C}$ and $\mathcal{U}$ a strong contravariant finite subcategory in $\mathcal{X}$.
The functor $\mathbb{F}: \mathcal{X}^{\mathcal{U}}_{L}/[\mathcal{X}]\longrightarrow \mathrm{mod}\mathcal{\underline{X}}$ is an equivalence of categories.
\ed{Th}

\Pf. (1) Firstly, we prove that the functor $\mathbb{F}$ is dense.

For any $F\in\mathrm{mod}\mathcal{\underline{X}}$, then there is an exact sequence in $\mathrm{mod}\mathcal{\underline{X}}$:
$$(\ddag1):\xymatrix{ \mathrm{Hom}_{\mathcal{\underline{X}}}(-,X_{1})\ar[r]^{u}
& \mathrm{Hom}_{\mathcal{\underline{X}}}(-,X_{0})\ar[r]&F\ar[r]&0}$$
with $X_{0}$, $X_{1}\in \mathcal{X}$. By Yoneda's Lemma, we have the follows isomorphism:
$$\mathrm{Hom}_{\mathcal{\underline{X}}}(\mathrm{Hom}_{\mathcal{\underline{X}}}(-,X_{1}),\mathrm{Hom}_{\mathcal{\underline{X}}}(-,X_{0}))\cong\mathrm{Hom}_{\mathcal{\underline{X}}}(X_{1},X_{0}).$$
Then there is a morphism $a$: $X_{1}\longrightarrow X_{0}$ such that $u=\mathrm{Hom}_{\mathcal{\underline{X}}}(-,a)$.
Since $\mathcal{U}$ is strong contravariant finite in $\mathcal{X}$, there is an $\mathbb{E}$-triangle
$$\xymatrix{A\ar[r] &U\ar[r]^{b}&X_{0}\ar@{.>}[r]&}$$
with $b$ right $\mathcal{U}$-approximation. By \cite[Corollary~3.16]{NP}, there exists an $\mathbb{E}$-triangle
$$(\ddag2):\xymatrix{B\ar[r] &X_{1}\bigoplus U\ar[r]^{~~(a,b)}&X_{0}\ar@{.>}[r]&.}$$
It is easy to see that the sequence $(\ddag2)$ is $\mathcal{U}$-exact since $b$ is right $\mathcal{U}$-approximation.
It follows from the sequence $(\ddag2)$ that $B\in\mathcal{X}^{\mathcal{U}}_{L}$.
By Lemma \ref{lemma1}, there is an exact sequence in $\mathrm{mod}\mathcal{\underline{X}}$
$$(\ddag3):\xymatrix{\mathrm{Hom}_{\mathcal{\underline{X}}}(-,X_{0})\ar[rr]^{\mathrm{Hom}_{\mathcal{\underline{X}}}(-,a)=u}&&\mathrm{Hom}_{\mathcal{\underline{X}}}(-,X_{1})\ar[r]&\mathbb{F}(B)\ar[r]&0.}$$
From the sequences $(\ddag1)$ and $(\ddag3)$, $\mathbb{F}(B)\cong F$.

(2) Secondly, we prove that the functor $\mathbb{F}$ is full.

Let $A$, $B\in \mathcal{X}^{\mathcal{U}}_{L}$ and $\alpha\in\mathrm{Hom}_{\mathrm{mod}\mathcal{\underline{X}}}(\mathbb{F}A,\mathbb{F}B)$.
By the definition of $\mathcal{X}^{\mathcal{U}}_{L}$, there are two $\mathcal{U}$-exact $\mathbb{E}$-triangles:
$$\xymatrix{A\ar[r]^{e} &X_{0}\ar[r]^{c}&X_{1}\ar@{.>}[r]&}$$
and
$$(\ddag4):\xymatrix{B\ar[r]^{g} &Y_{0}\ar[r]^{d}&Y_{1}\ar@{.>}[r]&}$$
with $X_{i}$, $Y_{i}\in \mathcal{X}$ for $i=0$, $1$.
Notice that $\mathrm{Hom}_{\mathcal{\underline{X}}}(-,X_{1})$ is projective in $\mathrm{mod}\mathcal{\underline{X}}$.
By Lemma \ref{lemma1}, we can get the following commutative diagram with exact rows in $\mathrm{mod}\mathcal{\underline{X}}$.
$$\xymatrix{\mathrm{Hom}_{\mathcal{\underline{X}}}(-,X_{0})\ar[r]\ar[d]^{v_{0}}&\mathrm{Hom}_{\mathcal{\underline{X}}}(-,X_{1})\ar[r]\ar[d]^{v_{1}}&\mathbb{F}(A)\ar[r]\ar[d]^{\alpha}&0\\
\mathrm{Hom}_{\mathcal{\underline{X}}}(-,Y_{0})\ar[r]&\mathrm{Hom}_{\mathcal{\underline{X}}}(-,Y_{1})\ar[r]&\mathbb{F}(B)\ar[r]&0
}$$
By Yoneda's Lemma, there are two morphisms $\underline{f_{i}}\in\mathrm{Hom}_{\mathcal{\underline{X}}}(X_{i},Y_{i})$ such that
$v_{i}=\mathrm{Hom}_{\mathcal{\underline{X}}}(-,\underline{f_{i}})$ and $\underline{f_{1}c}=\underline{df_{0}}$, where $i=0$, $1$.
Then there is an object $U\in \mathcal{U}$ and two morphisms $x\in\mathrm{Hom}_{\mathcal{C}}(U,Y_{1})$ and $y\in\mathrm{Hom}_{\mathcal{C}}(X_{0},U)$
such that the following diagram commutes:
$$\xymatrix{X_{0}\ar[dr]_{y}\ar[rr]^{f_{1}c-df_{0}}&&Y_{1}\\
&U\ar[ur]_{x}}$$
Then we have the following commutative diagram:
$$\xymatrix{
&&U\ar[d]^{x}\ar@{.>}[dl]_{z}&\\
B\ar[r]^{g} &Y_{0}\ar[r]^{d}&Y_{1}\ar@{.>}[r]&}$$
since the sequence $(\ddag4)$ is $\mathcal{U}$-exact. So we have that
$$f_{1}c-df_{0}=xy=dzy~~\mathrm{and}~~f_{1}c=dzy+df_{0}=d(zy+f_{0}).$$
Thus we can obtain the following $\mathbb{E}$-triangles commutative diagram by $(\mathrm{ET3})^{\mathrm{op}}$:
$$\xymatrix{
A\ar[r]^{e}\ar@{.>}[d]^{f}&X_{0}\ar[r]^{c}\ar[d]^{zy+f_{0}}&X_{1}\ar@{.>}[r]\ar[d]^{f_{1}}&\\
B\ar[r]^{g} &Y_{0}\ar[r]^{d}&Y_{1}\ar@{.>}[r]&}$$
Notice that the morphism $zy$ factors through $U\in\mathcal{U}$. Then $\underline{zy+f_{0}}=\underline{f_{0}}$.
Consequently, we have that $\alpha=\mathbb{F}([f])$. i.e., $\mathbb{F}$ is full.

(3) Finally, we prove that the functor $\mathbb{F}$ is faithful.

Let $h$: $A\longrightarrow B$ be in $\mathcal{X}^{\mathcal{U}}_{L}$ satisfying $\mathbb{F}([h])=0$. By Lemma \ref{lemma1},
we have that $0=\mathbb{F}([h])=\mathbb{H}([h])=\mathbb{E}(\mathcal{X},h)$.
Since $\mathcal{C}$ has enough injective objects, there is an $\mathbb{E}$-triangle
$$\xymatrix{A\ar[r]^{i} &I\ar[r]^{j}&A_{1}\ar@{.>}[r]&.}$$
By \cite[Corollary~3.16]{NP}, there exists an $\mathbb{E}$-triangle
$$(\ddag5):\xymatrix{A\ar[r]^{\binom{h}{i}~~~} &B\bigoplus I\ar[r]^{k}&A_{2}\ar@{.>}[r]&.}$$
Applying the functor $\mathrm{Hom}_{\mathcal{C}}(\mathcal{X},-)$ to the sequence $(\ddag5)$, we obtain a long exact sequence:
$$\xymatrix{\mathrm{Hom}_{\mathcal{C}}(\mathcal{X},B\bigoplus I)\ar[r]&\mathrm{Hom}_{\mathcal{C}}(\mathcal{X},A_{2})\ar[r]^{~~~\delta_{\sharp}}&
\mathbb{E}(\mathcal{X},A)\ar[rr]^{\mathbb{E}(\mathcal{X},\binom{h}{i})=0~~~}&&\mathbb{E}(\mathcal{X},B\bigoplus I).
}$$
Thus $\delta_{\sharp}$ is surjective. Since $A\in\mathcal{X}^{\mathcal{U}}_{L}$, there is a $\mathcal{U}$-exact $\mathbb{E}$-triangle
$$\xymatrix{A\ar[r]^{l} &X_{0}\ar[r]^{m}&X_{1}\ar@{.>}[r]&.}$$
with $X_{0}$, $X_{1}\in\mathcal{X}$.
Since $\delta_{\sharp}$ is surjective, there is a morphism $n\in \mathrm{Hom}_{\mathcal{C}}(X_{1},A_{2})$ such that
the following diagram is commutative:
$$\xymatrix{
A\ar[r]^{l}\ar@{=}[d]&X_{0}\ar[r]^{m}\ar[d]^{^{\binom{u}{v}~~~}}&X_{1}\ar@{.>}[r]\ar[d]^{n}&\\
A\ar[r]^{\binom{h}{i}~~~} &B\bigoplus I\ar[r]^{k}&A_{2}\ar@{.>}[r]&}$$
Thus $h=ul$ from the above diagram. i.e., $h$ factors through $X_{0}\in\mathcal{X}$.
Consequently, $[h]=0$. i.e., the functor $\mathbb{F}$ is faithful.
\ \hfill $\Box$

\vskip 10pt
In this section, if we take $\mathcal{U}=\mathcal{P}$, then we have the following result.

\bg{Cor}\cite[Theorem~3.4]{ZZ1}\label{Cor1}
Let $\mathcal{X}$ be a selforthogonal subcategory of an extriangulated category $\mathcal{C}$ with enough projective objects and injective objects.
The functor $\mathbb{F}: \mathcal{X}^{\mathcal{P}}_{L}/[\mathcal{X}]\longrightarrow \mathrm{mod}\mathcal{\underline{X}}$ is an equivalence of categories.
\ed{Cor}

Let $\mathcal{C}$ be an extriangulated category, $\mathcal{X}$ a subcategory of $\mathcal{C}$.
Then $\mathcal{X}$ is called cluster-tilting \cite{CZZ} if it satisfies the following conditions:

(1) $\mathcal{X}$ is strong functorially finite in $\mathcal{C}$;

(2) $M\in\mathcal{X}$ if and only if $\mathbb{E}(M,\mathcal{X})=0$;

(3) $M\in\mathcal{X}$ if and only if $\mathbb{E}(\mathcal{X},M)=0$.

\bg{Cor}
Let $\mathcal{C}$ be an extriangulated category with enough projective objects and injective objects,
and $\mathcal{X}$ a cluster-tilting subcategory of $\mathcal{C}$. Then $\mathcal{C}/[\mathcal{X}]\cong \mathrm{mod}\mathcal{\underline{X}}$.
\ed{Cor}

\Pf. By Remark 2.11 in \cite{ZZ1}, we have that $\mathcal{X}^{\mathcal{P}}_{L}=\mathcal{C}$. It follows from Corollary \ref{Cor1}.
\ \hfill $\Box$

\bg{Cor}
Let $\mathcal{X}$ be a strong contravariantly finite selforthogonal subcategory of $\mathcal{C}$ and $\mathcal{U}$ a strong contravariant finite subcategory in $\mathcal{X}$.
Then $\mathcal{X}^{\mathcal{U}}_{L}/[\mathcal{X}]$ is an abelian category.
\ed{Cor}

\Pf. It follows from Theorem \ref{Th1} and Proposition \ref{prop1}.
\ \hfill $\Box$

\section{Quotients of extriangulated categories II}
%

In this section, we mainly present the latter half of Theorem \ref{Th0}.
Let $\mathcal{C}$ be an extriangulated category and $\mathcal{X}$ a selforthogonal subcategory of $\mathcal{C}$.
Assume that $\mathcal{V}$ is a strong covariantly finite subcategory in $\mathcal{X}$.
The stable category $\mathcal{\overline{X}}$ of $\mathcal{X}$ is the quotient category $\mathcal{X}/[\mathcal{V}]$,
The category which has same objects as $\mathcal{X}$ and morphisms are defined as follows
$$\mathrm{Hom}_{\mathcal{\overline{X}}}(A,B)=\mathrm{Hom}_{\mathcal{X}}(A,B)/[\mathcal{V}](A,B),$$
and we denote by $\overline{f}$ the residue class in $\mathrm{Hom}_{\mathcal{\overline{X}}}(A,B)$ of a morphism $f$ of $\mathrm{Hom}_{\mathcal{X}}(A,B)$.

We denote by $\mathcal{X}^{R}_{\mathcal{V}}$ the subcategory of objects $A$ of $\mathcal{C}$ which admits an $\mathbb{E}$-triangle:
$$\xymatrix{X_{1}\ar[r] &X_{0}\ar[r]&A\ar@{.>}[r]&}$$ with $X_{0}$ and $X_{1}$ in $\mathcal{X}$
and it is $\mathcal{V}$-coexact. Since $\mathcal{V}$ is a strong covariantly finite subcategory in $\mathcal{X}$, for $X\in \mathcal{X}$,
there is an $\mathbb{E}$-triangle
$$\xymatrix{X\ar[r]^{x} &V\ar[r]&B\ar@{.>}[r]&}$$
with $V$ in $\mathcal{V}$ and $x$ left $\mathcal{V}$-approximation.
Denoted by $\Omega^{-1}_{\mathcal{V}}\mathcal{X}$ the subcategory of objects $B$ satisfies the above conditions.
It is easy to see that $\Omega^{-1}_{\mathcal{V}}\mathcal{X}\subseteq\mathcal{X}^{R}_{\mathcal{V}}$, thus we have the quotient category $\mathcal{X}^{R}_{\mathcal{V}}/[\Omega^{-1}_{\mathcal{V}}\mathcal{X}]$
and denote by $[f]$ the residue class in $\mathcal{X}^{R}_{\mathcal{V}}/[\Omega^{-1}_{\mathcal{V}}\mathcal{X}]$ of a morphism $f$ of $\mathcal{X}^{R}_{\mathcal{V}}$.

We consider the restriction of functor $\mathrm{Hom}(-, -)$ to $\mathcal{X}^{R}_{\mathcal{V}}$, denoted by $\mathbb{K}$:
$$\mathcal{X}^{R}_{\mathcal{V}}\longrightarrow \mathrm{Mod}\mathcal{\overline{X}}~;
A\longmapsto \mathbb{K}(A)=\mathrm{Hom}_{\mathcal{\overline{C}}}(-,A)|_{\mathcal{\overline{X}}}.$$
We claim that $\mathbb{K}(A)=0$ for any $A\in\Omega^{-1}_{\mathcal{V}}\mathcal{X}$.
Indeed, there is an $\mathbb{E}$-triangle
$$\xymatrix{X\ar[r]^{x} &V\ar[r]^{y}&A\ar@{.>}[r]&}$$
with $V$ in $\mathcal{V}$ and $x$ left $\mathcal{V}$-approximation since $A\in\Omega^{-1}_{\mathcal{V}}\mathcal{X}$.
Applying the functor $\mathrm{Hom}_{\mathcal{C}}(\mathcal{X},-)$ to above $\mathbb{E}$-triangle, we have the following long exact sequence
$$\xymatrix{\mathrm{Hom}_{\mathcal{C}}(\mathcal{X},X)\ar[r] &\mathrm{Hom}_{\mathcal{C}}(\mathcal{X},V)\ar[rr]^{\mathrm{Hom}_{\mathcal{C}}(\mathcal{X},y)}&&
\mathrm{Hom}_{\mathcal{C}}(\mathcal{X},A)\ar[r]&\mathbb{E}(\mathcal{X},X)=0.}$$
Thus $\mathrm{Hom}_{\mathcal{C}}(\mathcal{X},y)$ is surjective since $\mathcal{X}$ is selforthogonal, i.e., for any $f \in\mathrm{Hom}_{\mathcal{C}}(\mathcal{X},A)$,
$f$ factors through $y$, and then $\overline{f}=0.$ Consequently, $\mathbb{K}(A)=0$ for any $A\in\Omega^{-1}_{\mathcal{V}}\mathcal{X}$.

Let $\pi'$: $\mathcal{X}^{R}_{\mathcal{V}}\longrightarrow \mathcal{X}^{R}_{\mathcal{V}}/[\Omega^{-1}_{\mathcal{V}}\mathcal{X}]$ be projective functor.
Since $\mathbb{K}(A)=0$ for any $A\in\Omega^{-1}_{\mathcal{V}}\mathcal{X}$. By the universal property of $\pi'$, there is a functor
$$\mathbb{G}: \mathcal{X}^{R}_{\mathcal{V}}/[\Omega^{-1}_{\mathcal{V}}\mathcal{X}]\longrightarrow \mathrm{Mod}\mathcal{\overline{X}}$$
such that the following diagram commutes:
$$\xymatrix{
\mathcal{X}^{R}_{\mathcal{V}}\ar[rr]^{\pi'}\ar[dr]_{\mathbb{K}}&&\mathcal{X}^{R}_{\mathcal{V}}/[\Omega^{-1}_{\mathcal{V}}\mathcal{X}]\ar@{.>}[dl]^{\mathbb{G}}\\
&\mathrm{Mod}\mathcal{\overline{X}}&.
}$$

\vskip 10pt
In order to complete the proof of Theorem \ref{Th0}, the following lemma is important.

\bg{Lem}\label{lemma2}
Let $\mathcal{X}$ be a selforthogonal subcategory of an extriangulated category $\mathcal{C}$ and $\mathcal{V}$ a strong covariantly finite subcategory in $\mathcal{X}$.
For any $\mathcal{V}$-coexact $\mathbb{E}$-triangle
$$(\natural):\xymatrix{X_{1}\ar[r]^{a} &X_{0}\ar[r]^{b}&A\ar@{.>}[r]&}$$
with $X_{0}$, $X_{1}\in \mathcal{X}$, then there is an exact sequence in $\mathrm{Mod}\mathcal{\overline{X}}$
$$\xymatrix{\mathrm{Hom}_{\mathcal{\overline{C}}}(-,X_{1})\ar[r]^{\mathbb{K}(a)}&\mathrm{Hom}_{\mathcal{\overline{C}}}(-,X_{0})\ar[r]^{\mathbb{K}(b)}
&\mathrm{Hom}_{\mathcal{\overline{C}}}(-,A)\ar[r]&0.}$$
Thus we have that $\mathbb{K}(A)=\mathbb{G}(A)\in \mathrm{mod}\mathcal{\underline{X}}$.
\ed{Lem}

\Pf. For any $X\in \mathcal{X}$, applying the functor $\mathrm{Hom}_{\mathcal{C}}(X,-)$ to above $\mathbb{E}$-triangle $(\natural)$, we have the following long exact sequence
$$\xymatrix{\mathrm{Hom}_{\mathcal{C}}(X,X_{1})\ar[r] &\mathrm{Hom}_{\mathcal{C}}(X,X_{0})\ar[r]&
\mathrm{Hom}_{\mathcal{C}}(X,A)\ar[r]&\mathbb{E}(X,X_{1})=0}$$
since $\mathcal{X}$ is selforthogonal. It is easy to prove that $b$ is a right $\mathcal{X}$-approximation and $\mathbb{K}(b)$ is surjective.

For any $c\in \mathrm{Hom}_{\mathcal{C}}(X,X_{0})$ satisfying $\mathbb{K}(c)=\overline{bc}=0$ in $\overline{\mathcal{X}}$, then we can get the following
commutative diagram for some $V\in \mathcal{V}$:
$$\xymatrix{
&X\ar[r]^{d}\ar[d]^{c}&V\ar[d]^{e}&\\
X_{1}\ar[r]^{a} &X_{0}\ar[r]^{b}&A\ar@{.>}[r]&
}$$
Since $b$ is a right $\mathcal{X}$-approximation, then we have the following commutative diagram:
$$\xymatrix{
&&V\ar@{.>}[dl]_{f}\ar[d]^{e}&\\
X_{1}\ar[r]^{a} &X_{0}\ar[r]^{b}&A\ar@{.>}[r]&
}$$
Thus $bfd=ed=bc$ and $b(fd-c)=0$. Then we have the following commutative diagram:
$$\xymatrix{
X\ar@{.>}[d]^{g}\ar[dr]^{fd-c}&&&\\
X_{1}\ar[r]^{a} &X_{0}\ar[r]^{b}&A\ar@{.>}[r]&
}$$
So $fd-c=ag$ and $c=ag+fd$. Notice that $\overline{fd}=0$. Then $\overline{c}=\overline{ag}$.
i.e., $\mathrm{Ker}\mathbb{K}(b)=\mathrm{Im}\mathbb{K}(a)$.
Consequently, we can obtain the following exact sequence in $\mathrm{Mod}\mathcal{\overline{X}}$
$$\xymatrix{\mathrm{Hom}_{\mathcal{\overline{C}}}(-,X_{1})\ar[r]^{\mathbb{K}(a)}&\mathrm{Hom}_{\mathcal{\overline{C}}}(-,X_{0})\ar[r]^{\mathbb{K}(b)}
&\mathrm{Hom}_{\mathcal{\overline{C}}}(-,A)\ar[r]&0.}$$
\ \hfill $\Box$

Finally, we give the proof of the latter half of Theorem 1.1.

\bg{Th}\label{Th2}
Let $\mathcal{X}$ be a selforthogonal subcategory of an extriangulated category $\mathcal{C}$ and $\mathcal{V}$ a strong covariantly finite subcategory in $\mathcal{X}$.
The functor $\mathbb{G}: \mathcal{X}^{R}_{\mathcal{V}}/[\Omega^{-1}_{\mathcal{V}}\mathcal{X}]\longrightarrow \mathrm{mod}\mathcal{\overline{X}}$ is an equivalence of categories.
\ed{Th}

\Pf. (1) Firstly, we prove that the functor $\mathbb{G}$ is dense.

For any $F\in \mathrm{mod}\mathcal{\overline{X}}$, then there is an exact sequence
$$(\natural 1): \xymatrix{\mathrm{Hom}_{\mathcal{\overline{C}}}(-,X_{1})\ar[r]^{u}&\mathrm{Hom}_{\mathcal{\overline{C}}}(-,X_{0})\ar[r]
&F\ar[r]&0.}$$
with $X_{0}$, $X_{1}\in \mathcal{X}$. By Yoneda's Lemma, there is a morphism $a$: $X_{1}\longrightarrow X_{0}$ such that
$u=\mathbb{K}(a)$. Since $\mathcal{V}$ is a strong covariantly finite subcategory in $\mathcal{X}$,
there is an $\mathbb{E}$-triangle
$$\xymatrix{X_{1}\ar[r]^{b} &V_{1}\ar[r]&M\ar@{.>}[r]&}$$
with $b$ a left $\mathcal{V}$-approximation.
By \cite[Corollary~3.16]{NP}, we have the following $\mathbb{E}$-triangle
$$(\natural 2): \xymatrix{X_{1}\ar[r]^{\binom{a}{b}~~~} &X_{0}\bigoplus V_{1}\ar[r]&M\ar@{.>}[r]&.}$$
Since $b$ is a left $\mathcal{V}$-approximation, it is easy to prove that the $\mathbb{E}$-triangle $(\natural 2)$ is $\mathcal{V}$-coexact,
and then $N\in\mathcal{X}^{R}_{\mathcal{V}}$ by the definition of $\mathcal{X}^{R}_{\mathcal{V}}$.
By Lemma \ref{lemma2}, we have the following exact sequence
$$(\natural 3): \xymatrix{\mathrm{Hom}_{\mathcal{\overline{C}}}(-,X_{1})\ar[r]^{u~~~~~~~~~~~~~~~}&\mathrm{Hom}_{\mathcal{\overline{C}}}(-,X_{0}\bigoplus V_{1})\cong \mathrm{Hom}_{\mathcal{\overline{C}}}(-,X_{0})\ar[r]
&\mathrm{Hom}_{\mathcal{\overline{C}}}(-,N)\ar[r]&0}$$
since $X_{0}\bigoplus V_{1}\cong X_{0}$ in $\overline{\mathcal{C}}$.
From the sequences $(\natural 1)$ and $(\natural 3)$, we have that $F\cong \mathrm{Hom}_{\mathcal{\overline{C}}}(-,N)=\mathbb{K}(N)=\mathbb{G}(N)$.
i.e., the functor $\mathbb{G}$ is dense.

(2) Secondly, we prove that the functor $\mathbb{G}$ is full.

For any $A$, $B\in \mathcal{X}^{R}_{\mathcal{V}}$ and $\alpha\in \mathrm{Hom}_{\mathrm{mod}\mathcal{\overline{X}}}(\mathbb{K}(A),\mathbb{K}(B))$.
Then we have the following $\mathbb{E}$-triangles:
$$(\natural 4): \xymatrix{X_{1}\ar[r]^{c} &X_{0}\ar[r]^{d}&A\ar@{.>}[r]&}$$
and
$$\xymatrix{Y_{1}\ar[r]^{e} &Y_{0}\ar[r]^{f}&B\ar@{.>}[r]&}$$
with $X_{i}$, $Y_{i}\in \mathcal{X}$ for $i=0$, $1$ and the two $\mathbb{E}$-triangles are $\mathcal{V}$-coexact.
Notice that $\mathbb{K}(X_{0})$ is projective in $\mathrm{mod}\mathcal{\overline{X}}$.
By Lemma \ref{lemma2}, we can get the following commutative diagram with exact rows in $\mathrm{mod}\mathcal{\overline{X}}$.
$$\xymatrix{\mathbb{K}(X_{1})\ar[r]^{\mathbb{K}(c)}\ar[d]^{\gamma}&\mathbb{K}(X_{0})\ar[r]^{\mathbb{K}(d)}\ar[d]^{\beta}&\mathbb{K}(A)\ar[r]\ar[d]^{\alpha}&0\\
\mathbb{K}(Y_{1})\ar[r]^{\mathbb{K}(e)}&\mathbb{K}(Y_{0})\ar[r]^{\mathbb{K}(f)}&\mathbb{K}(B)\ar[r]&0
}$$
By Yoneda's Lemma, there are two morphisms $h\in\mathrm{Hom}_{\mathcal{C}}(X_{1},Y_{1})$ and $g\in\mathrm{Hom}_{\mathcal{C}}(X_{0},Y_{0})$ such that
$\gamma=\mathbb{K}(h)$, $\beta=\mathbb{K}(g)$ and $\overline{gc}=\overline{eh}$.
Then there is an object $V\in \mathcal{V}$ and two morphisms $m\in\mathrm{Hom}_{\mathcal{C}}(V,Y_{0})$ and $n\in\mathrm{Hom}_{\mathcal{C}}(X_{1},V)$
such that the following diagram commutates:
$$\xymatrix{X_{1}\ar[dr]_{n}\ar[rr]^{gc-eh}&&Y_{0}\\
&V\ar[ur]_{m}}$$
Since the sequence $(\natural 4)$ is $\mathcal{V}$-coexact, we have the following commutative diagram:
$$\xymatrix{
V&&&\\
X_{1}\ar[u]^{n}\ar[r]^{c} &X_{0}\ar@{.>}[ul]_{i}\ar[r]^{d}&A\ar@{.>}[r]&}$$
So we have that
$$gc-eh=mn=mic~~\mathrm{and}~~eh=gc-mic=(g-mi)c.$$
Thus we can obtain the following $\mathbb{E}$-triangles commutative diagram by $(\mathrm{ET3})$:
$$\xymatrix{
X_{1}\ar[r]^{c}\ar[d]^{h}&X_{0}\ar[r]^{d}\ar[d]^{g-mi}&A\ar@{.>}[r]\ar@{.>}[d]^{j}&\\
Y_{1}\ar[r]^{e} &Y_{0}\ar[r]^{f}&B\ar@{.>}[r]&}$$
Notice that the morphism $mi$ factors through $V\in\mathcal{V}$. Then $\overline{g-mi}=\overline{g}$.
Consequently, we have that $\alpha=\mathbb{G}([j])$. i.e., the functor $\mathbb{G}$ is full.

(3) Finally, we prove that the functor $\mathbb{G}$ is faithful.

Let $x$: $A\longrightarrow B$ be in $\mathcal{X}^{R}_{\mathcal{V}}$ satisfying $\mathbb{G}([x])=0$.
Then $0=\mathbb{G}([x])=\mathbb{K}([x])=\mathrm{Hom}_{\overline{\mathcal{C}}}(\mathcal{\overline{X}},x)$.
Since $A\in\mathcal{X}^{R}_{\mathcal{V}}$ and $\mathcal{V}$ is a strong covariantly finite subcategory in $\mathcal{X}$, there
are two $\mathbb{E}$-triangles
$$(\natural 5): \xymatrix{X_{1}\ar[r]^{y} &X_{0}\ar[r]^{z}&A\ar@{.>}[r]&}$$
and
$$\xymatrix{X_{0}\ar[r]^{u} &V_{0}\ar[r]^{v}&X'_{0}\ar@{.>}[r]&}$$
where $X_{0}$, $X_{1}\in X$, $u$ is left $\mathcal{V}$-approximation and the $\mathbb{E}$-triangle $(\natural 5)$ is $\mathcal{V}$-coexact.
Then we have the following diagram is commutative by $(\mathrm{ET4})$:
$$\xymatrix{
X_{1}\ar@{=}[d]\ar[r]^{y} &X_{0}\ar[d]^{u}\ar[r]^{z}&A\ar@{.>}[r]\ar[d]^{i}&\\
X_{1}\ar[r]^{t}&V_{0}\ar[r]^{j}\ar[d]^{v}&J\ar@{.>}[r]\ar[d]^{k}&\\
&X_{0}'\ar@{.>}[d]\ar@{=}[r]&X_{0}'\ar@{.>}[d]\\
&&&
}$$
Since $\mathrm{Hom}_{\overline{\mathcal{C}}}(\mathcal{\overline{X}},x)=0$, $\overline{xz}=0$.
Then we have the following commutative diagram:
$$\xymatrix{
X_{1}\ar[r]^{y} &X_{0}\ar[d]^{n}\ar[r]^{z}&A\ar@{.>}[r]\ar[d]^{x}&\\
&V_{1}\ar[r]^{m}&B&
}$$
where $V_{1}\in \mathcal{V}$.
Since $u$ is left $\mathcal{V}$-approximation, we have the following commutative diagram:
$$\xymatrix{
V_{1}&&&\\
X_{0}\ar[u]^{n}\ar[r]^{u} &V_{0}\ar@{.>}[ul]_{l}\ar[r]^{v}&X'_{0}\ar@{.>}[r]&
}$$
So $xz=mn=mlu$. By \cite[Lemma~ 3.13]{NP}, there is a morphism $q$: $J\longrightarrow B$ which makes
the following diagram commutative:
$$\xymatrix{
X_{0}\ar[d]^{u}\ar[r]^{z}&A\ar[ddr]^{x}\ar[d]^{i}&\\
V_{0}\ar[r]^{j}\ar[drr]_{ml}&J\ar@{.>}[dr]^{q}&\\
&&B
}$$
Thus $x=qi$. i.e., the morphism $x$ factors through $J$. Next we will prove that $J\in \Omega^{-1}_{\mathcal{V}}\mathcal{X}$,
and then $[x]=0$. We consider the $\mathbb{E}$-triangle
$$\xymatrix{X_{1}\ar[r]^{t} &V_{0}\ar[r]^{j}&J\ar@{.>}[r]&.}$$
We only need to prove that $t$ is a left $\mathcal{V}$-approximation.
Since the $\mathbb{E}$-triangle $(\natural 5)$ is $\mathcal{V}$-coexact and $u$ is left $\mathcal{V}$-approximation,
we have the following commutative diagram for any $V\in \mathcal{V}$ and $p$: $X_{1}\longrightarrow V$:
$$\xymatrix{
X_{1}\ar@{=}[d]\ar[r]^{y}&X_{0}\ar@{.>}[ddr]^{s}\ar[d]^{u}&\\
X_{1}\ar[r]^{t}\ar[drr]_{p}&V_{0}\ar@{.>}[dr]^{v}&\\
&&V
}$$
Thus $p=sy=vuy=vt$, and then $t$ is a left $\mathcal{V}$-approximation.
Consequently, $[x]=0$. i.e., the functor $\mathbb{G}$ is faithful.
\ \hfill $\Box$

\bg{Cor}\label{}
$\Omega^{-1}_{\mathcal{V}}\mathcal{X}=\{A\in \mathcal{X}^{R}_{\mathcal{V}}~|~ \mathrm{Hom}_{\mathcal{\overline{C}}}(\mathcal{X},A)=0\}.$
\ed{Cor}

\Pf. If $A\in\Omega^{-1}_{\mathcal{V}}\mathcal{X}$, then there is an $\mathbb{E}$-triangle
$$\xymatrix{X\ar[r]^{x} &V\ar[r]&A\ar@{.>}[r]&}$$
with $V\in\mathcal{V}$, $X\in \mathcal{X}$ and $x$ left $\mathcal{V}$-approximation.
Applying the functor $\mathrm{Hom}_{\mathcal{C}}(\mathcal{X},-)$ to the above $\mathbb{E}$-triangle, we can get the following exact sequence
$$\xymatrix{\mathrm{Hom}_{\mathcal{C}}(\mathcal{X},V)\ar[r] &\mathrm{Hom}_{\mathcal{C}}(\mathcal{X},A)\ar[r]&\mathbb{E}(\mathcal{X},X)=0.}$$
Thus for any $f\in \mathrm{Hom}_{\mathcal{C}}(X,A)$ with $X\in\mathcal{X}$ factors through $V\in\mathcal{V}$, i.e., $\mathrm{Hom}_{\mathcal{\overline{C}}}(\mathcal{X},A)=0$.

Conversely, if $A\in\mathcal{X}^{R}_{\mathcal{V}}$ and $\mathrm{Hom}_{\mathcal{\overline{C}}}(\mathcal{X},A)=0$,
then $0=\mathrm{Hom}_{\mathcal{\overline{C}}}(\mathcal{X},A)=\mathbb{G}(A)$.
So $\mathbb{G}(1_{A})=0$, and then $1_{A}=0$ in $\mathcal{X}^{R}_{\mathcal{V}}/\Omega^{-1}_{\mathcal{V}}\mathcal{X}$ since $\mathbb{G}$ is faithful.
Thus $A\in\Omega^{-1}_{\mathcal{V}}\mathcal{X}$.
\ \hfill $\Box$

\vskip 10pt
In this section, if we take $\mathcal{V}=\mathcal{I}$, then we have the following result.

\bg{Cor}\cite[Theorem~3.13]{ZZ1}\label{}
Let $\mathcal{X}$ be a selforthogonal subcategory of an extriangulated category $\mathcal{C}$ with enough projective objects and injective objects.
The functor $\mathbb{G}: \mathcal{X}^{R}_{\mathcal{I}}/[\Omega^{-1}_{\mathcal{I}}\mathcal{X}]\longrightarrow \mathrm{mod}\mathcal{\overline{X}}$ is an equivalence of categories.
\ed{Cor}

\section{Declarations}

\textbf{Conflict of interest} This work does not have any conflicts of interest.


{\small

}


\begin{thebibliography}{17}

\bibitem{Au}  M. Auslander. Coherent functors. In Proceedings Conference Categorical Algebra
(La Jolla). Springer-Verlag, Berlin-Heidelberg-NewYork, 1966, 189--231.

\bibitem{Au1} M. Auslander. Representation Dimension of Artin Algebras. Queen Mary College Mathematics Notes. London: Queen Mary College, 1971.

\bibitem{BMR} B. Buan, R. Marsh, and I. Reiten. Cluster-tilted algebras. Trans. Amer. Math.Soc., 2007, 359(1): 323--332.

\bibitem{CZZ} W. Chang, P. Zhou and B. Zhu. Cluster subalgebras and cotorsion pairs in Frobenius extriangulated categories. Algebr. and Rep. Theory, 2019, 22(5): 1051--1081.

\bibitem{DL}  L. Demonet and Y. Liu. Quotients of exact categories by cluster tilting subcategories as module categories. J. Pure Appl. Algebra, 2013, 217(12): 2282--2297.

\bibitem{HZ} Y. Hu and P. Zhou, Recollements arising from cotorsion pairs on extriangulated categories. Front. Math. China, 2021, 16(4): 937--955.

\bibitem{HZZ1} J. Huang, D. Zhang and P. Zhou, Two new classes of n-exangulated categories, Journal of Algebra, 2021, 568: 1--21.

\bibitem{HZZ} J. Huang, D. Zhang and P. Zhou. Proper classes and Goensteinness in extriangulated categories. Journal of Algebra, 2020, 551: 23--60.


\bibitem{LN} Y. Liu and H. Nakaoka, Hearts of twin cotorsion pairs on extriangulated categories, Journal of Algebra, 2019, 528: 96--149.

\bibitem{NP} H. Nakaoka and Y. Palu. Extriangulated categories, Hovey twin cotorsion pairs and model structures, Cahiers de Topologie et Geometrie Differentielle Categoriques. 2019, Volume LX-2, 117--193.

\bibitem{WWZ} L. Wang, J. Wei and H. Zhang, Recolloments of extriangulated categories. Colloquium Mathematicum, 2022, 167(2): 239--259.

\bibitem{ZB} B. Zhu and X. Zhuang, Tilting subcategories in extriangulated categories, Front. Math. China, 2020, 15(1), 225--253.

\bibitem{ZZ}  P. Zhou and B. Zhu. Triangulated quotient categories revisited. J.Algebra, 2018, 50: 196--232.

\bibitem{ZZ1}  P. Zhou and B. Zhu. Cluster-tilting subcategories in extriangulated categories. Theory and Applications of Categories, 2019, 34(8): 221--242.



\end{thebibliography}
\end{document}